\documentclass[twoside,a4paper,reqno,11pt]{amsart}
\usepackage{amsfonts, amsbsy, amsmath, amssymb, latexsym}
\usepackage{mathrsfs,array}
\usepackage[top=30mm,right=30mm,bottom=30mm,left=30mm]{geometry}
\usepackage{stmaryrd}
\usepackage{bm}
\usepackage{xcolor}
\usepackage{hyperref}
\usepackage{tikz-cd}
\usepackage{relsize}

\makeatletter
\@namedef{subjclassname@2020}{\textup{2020} Mathematics Subject Classification}
\makeatother

\def\thefootnote{\fnsymbol{footnote}}

 \renewcommand{\to}{\rightarrow}

\DeclareMathOperator{\Lie}{Lie}

\DeclareMathOperator{\sep}{s}

\DeclareMathOperator{\Rad}{Rad}

\makeatletter
\newcommand{\imod}[1]{\allowbreak\mkern4mu({\operator@font mod}\,\,#1)}
\makeatother

\newtheorem{theorem}{Theorem}[section]
\newtheorem{lemma}[theorem]{Lemma}
\newtheorem{proposition}[theorem]{Proposition}

\newtheorem{corollary}[theorem]{Corollary}

\theoremstyle{definition}

\newtheorem{example}[theorem]{Example}
\newtheorem{remark}[theorem]{Remark}

\begin{document}

\title{Maximal toroids and Cartan subgroups of algebraic groups}
\author[D.\ Sercombe]{Damian Sercombe}
\address
{Mathematisches Institut, Universität Freiburg, Ernst-Zermelo-Straße 1, 79104 Freiburg, Germany}
\email{damian.sercombe@math.uni-freiburg.de}

\begin{abstract} We introduce a unified theory of Cartan subgroups and maximal toroids -- defined as connected multiplicative type subgroups that are maximal amongst all such subgroups -- which holds for all affine algebraic groups over a field, regardless of smoothness. For instance we show that maximal toroids always exist, that they are invariant under base change, and that they are in natural 1--1 correspondence with Cartan subgroups. Our results generalise known results for Cartan subgroups and maximal tori of smooth affine algebraic groups, as well as their analogues for restricted Lie algebras. We conclude with some applications to, and a brief discussion of, some generation problems for algebraic groups.
\end{abstract}

\let\thefootnote\relax

\maketitle

\section{Introduction}

\noindent The structure theory of smooth algebraic groups over a field is well-established. A central role is played by the existence, conjugacy, and invariance under base change, of maximal tori. 
This allows one to rigidify a smooth algebraic group $G$ by fixing a maximal torus $T$ (or, equivalently, a Cartan subgroup), and then to study associated structures of $G$ relative to $T$ in a manner independent of the choice of $T$. Examples of such structures include the root system, root groups, the Weyl group, Borel and parabolic subgroups, etc. 
This rigidification is an important tool; for instance, it is a key step in the classification of reductive groups.

\vspace{2mm}\noindent The theory of finite-dimensional restricted Lie algebras over a field is more complicated, but proceeds along similar lines. Any such Lie algebra $\mathfrak{g}$ admits at least one Cartan subalgebra. Over an algebraically closed field of characteristic zero, any two Cartan subalgebras of $\mathfrak{g}$ are conjugate by an inner automorphism of $\mathfrak{g}$. 
In positive characteristic, the analogous statement fails. Indeed, it is not even true that all Cartan subalgebras of $\mathfrak{g}$ must have the same dimension. One circumvents this difficulty by fixing a Cartan subalgebra of maximal toral dimension in $\mathfrak{g}$; 
it is a result of Premet \cite[Thm.\ 1]{Pr} (building on work of Winter \cite{Wi}) that these are all equivalent by means of a finite number of operations known as ``elementary switchings". 
This rigidification of $\mathfrak{g}$ plays a similar role as with smooth algebraic groups, and is used in the classification of simple Lie algebras. 

\vspace{2mm}\noindent In this paper we introduce and develop a unified theory of Cartan subgroups and maximal toroids -- a generalisation of (maximal) tori/toral subalgebras -- which are applicable to any \textbf{affine} algebraic group, regardless of smoothness. Our definitions and results generalise those for smooth affine algebraic groups, and for restricted Lie algebras (via the Demazure-Gabriel correspondence with height 1 algebraic groups).

\vspace{2mm}\noindent Many results which hold in both the smooth and height 1 cases generalise remarkably well to an arbitrary affine algebraic group, giving us some reassurance that our definitions have been chosen correctly. For instance, in Theorem \ref{Cartans} we obtain the expected characterisation of Cartan subgroups as centralisers of maximal toroids. One notable exception, illustrated in Example \ref{notgenbyCartans}, is that an affine algebraic group is not necessarily generated by its Cartan subgroups.

\vspace{2mm}\noindent Having established a theory of maximal toroids of an affine algebraic group $G$, one wants to find a distinguished subclass of maximal toroids 
which are large enough and also similar enough to each other -- conjugacy is too much to hope for -- that they can be used to appropriately rigidify $G$. Our leading candidate is those toroids $T$ of $G$ such that, for some large integer $r$, their $r$'th Frobenius kernel $T_r$ has maximal order amongst all toroids of $G_r$. Exploring this further is a subject for future work.

\vspace{2mm}\noindent This paper was originally motivated by certain generation problems for algebraic groups: i.e. given a particular set of subgroups of an algebraic group $G$, does it generate all of $G$? Of particular interest is the set of all nilpotent subgroups, and the set of all multiplicative type subgroups. Whilst exploring such problems it became clear the importance of constructing analogues to Cartan subgroups and maximal tori without any smoothness hypotheses, which developed into this paper. As such we conclude this paper with some applications to, and a brief discussion of, such generation problems. These problems will be explored further in a subsequent paper.

\vspace{2mm}\noindent Our setup henceforth is as follows. Let $k$ be a field, say of characteristic $p\geq 0$. Let $G$ be an affine algebraic $k$-group, that is, an affine group scheme of finite type over $k$. 
A \textit{subgroup} of $G$ will refer to a locally closed $k$-subgroup scheme. 

\vspace{2mm}\noindent Let us recall some basic definitions. A \textit{torus} is a smooth connected affine $k$-group of multiplicative type. A \textit{torus of $G$} is a subgroup of $G$ that is isomorphic to a torus, and a \textit{maximal torus of $G$} is a torus of $G$ which is not strictly contained in any other torus of $G$. 
We define analogues of these notions without any smoothness hypotheses, as follows.

\vspace{2mm}\noindent A \textit{toroid} is a connected affine $k$-group of multiplicative type. [N.B. Not to be confused with the notion of a ``toroidal algebraic group" introduced by Rosenlicht in \cite{Ro}.] For example, if $p>0$ then the $p$'th roots of unity $\mu_p$ is a toroid but it is not smooth. A \textit{toroid of $G$} is a subgroup of $G$ that is isomorphic to a toroid, and a \textit{maximal toroid of $G$} is a toroid of $G$ which is not strictly contained in any other toroid of $G$.

\vspace{2mm}\noindent The existence of a maximal torus of $G$ is obvious by dimension considerations. 
The existence of a maximal toroid is a more difficult question, which we answer in the affirmative:

\begin{theorem}\label{maxtoroids} Let $G$ be an affine algebraic $k$-group.

\vspace{1mm}\noindent (a) Every toroid of $G$ is contained in a maximal toroid of $G$. Hence $G$ contains at least one maximal toroid. 

\vspace{1mm}\noindent (b) Suppose $G$ is smooth. Then every maximal toroid of $G$ is smooth. In other words, the maximal toroids of $G$ are precisely the maximal tori of $G$. 
\end{theorem}

\noindent Note that Theorem \ref{maxtoroids} fails without the assumption that $G$ is affine. For example, suppose $p>0$ and $G$ is an ordinary elliptic curve. Then, for each integer $i\geq 1$, the $i$'th Frobenius kernel $G_i$ of $G$ is a toroid of order $p^i$, 
yet every affine subgroup of $G$ is finite. 

\vspace{2mm}\noindent A subgroup $C$ of $G$ is called a \textit{Cartan subgroup} if it is connected, nilpotent, and satisfies $C=N_G(C)^{\circ}$. 
[N.B. In Corollary \ref{Cartanscor3} we will show that if $G$ is smooth then its Cartan subgroups are automatically smooth. So our definition of a Cartan subgroup is indeed a generalisation of the usual definition for smooth $G$.]

\vspace{2mm}\noindent It is a classical result, see for instance \cite[11.13, Thm.\ 12.6]{Bo}, that for smooth $G$ there is a natural 1-1 correspondence between the set of maximal tori of $G$ and the set of (smooth) Cartan subgroups of $G$. In Corollary \ref{Cartanscor1} we show that this 1-1 correspondence generalises nicely to the not-necessarily-smooth situation. But first, we characterise Cartan subgroups of $G$. Let $\overline{k}$ be an algebraic closure of $k$.

\begin{theorem}\label{Cartans} Let $G$ be an affine algebraic $k$-group. Let $C$ be a subgroup of $G$. The following are equivalent:

\vspace{1mm}\noindent (i) $C$ is a Cartan subgroup of $G$.

\vspace{1mm}\noindent (ii) $C_{\overline{k}}$ is a Cartan subgroup of $G_{\overline{k}}$.

\vspace{1mm}\noindent (iii) $C$ is connected, and 

\vspace{-1mm}\begin{itemize} \item if $p=0$ the Lie algebra $\Lie(C)$ is a Cartan subalgebra of $\Lie(G)$, \item if $p>0$ there exists a positive integer $r$ such that, for every $i\geq r$, the $i$'th Frobenius kernel $C_i$ is a Cartan subgroup of $G_i$. 
\end{itemize} 

\vspace{1mm}\noindent (iv) $C=Z_G(T)^{\circ}$, for some maximal toroid $T$ of $G$.

\vspace{1mm}\noindent (v) $C$ is a connected nilpotent subgroup of $G$ which is maximal amongst all such subgroups of $G$, and every toroid of $G$ which normalises $C$ must centralise $C$.
\end{theorem} 

\noindent We now present some consequences of Theorems \ref{maxtoroids} and \ref{Cartans}. Given an affine algebraic $k$-group $H$, we let $Z(H)^{\circ}_s$ denote the (unique) maximal central toroid of $H$.

\begin{corollary}\label{Cartanscor1} Let $G$ be an affine algebraic $k$-group. The maps $T \mapsto Z_G(T)^{\circ}$ and $C \mapsto Z(C)^{\circ}_s$ are inverse bijections between the set of maximal toroids of $G$ and the set of Cartan subgroups of $G$. 
\end{corollary}

\begin{corollary}\label{Cartanscor2} Every affine algebraic $k$-group contains at least one Cartan subgroup.
\end{corollary}

\begin{corollary}\label{Cartanscor3} Every Cartan subgroup of a smooth affine $k$-group is smooth.
\end{corollary}

\noindent It is an important result of Grothendieck \cite[XIV, Thm.\ 1.1]{SGA3} that a maximal torus of a smooth affine $k$-group $G$ remains a maximal torus after base change by $\smash{\overline{k}}/k$. Without the smoothness hypothesis on $G$, this statement fails. 
Quite unexpectedly, we find that the analogue of this statement for maximal toroids holds without any smoothness hypotheses on $G$ (in fact it generalises the aforementioned result, by Theorem \ref{maxtoroids}(b)). 

\begin{corollary}\label{Cartanscor4} Let $G$ be an affine algebraic $k$-group. Let $T$ be a subgroup of $G$. Then $T$ is a maximal toroid of $G$ if and only if $T_{\overline{k}}$ is a maximal toroid of $G_{\overline{k}}$.
\end{corollary}

\begin{corollary}\label{Cartanscor5} Assume $p>0$. Let $G$ be an affine algebraic $k$-group. Let $T$ be a subgroup of $G$. Then $T$ is a maximal toroid of $G$ if and only if $T$ is connected and there exists a positive integer $r$ such that, for every $i\geq r$, the $i$'th Frobenius kernel $T_i$ is a maximal toroid of $G_i$. 
\end{corollary}

\noindent We now turn our attention to some generation problems. It is well-known that any smooth connected affine $k$-group is generated by its Cartan subgroups, see for instance \cite[Thm.\ 13.3.6(ii), 13.3.7]{Sp}. 
The analogous statement fails without the smoothness assumption -- even over an algebraically closed field -- as illustrated in Example \ref{notgenbyCartans}. As such, we are interested in the subgroup $G_c$ of $G$ which is generated by its Cartan subgroups. 

\begin{theorem}\label{Gcsubgroup} Assume $p>0$. Let $G$ be an affine algebraic $k$-group. Let $G_c$ be the subgroup of $G$ generated by all Cartan subgroups of $G$. Then $G_c$ contains the Frobenius kernel $G_1$ of $G$. 
\end{theorem}

\noindent We would like to better understand this subgroup $G_c$ of $G$. For instance, is it true that $G=G_c$ if all unipotent subgroups of $G$ have height at most 1? Another question is: is $G_c$ normal in $G$? We suspect the answer to both questions is yes, but currently do not know.

\vspace{2mm}\noindent The rigidity of multiplicative type $k$-groups \cite[Thm.\ 12.36]{Mi} implies that any nilpotent connected affine algebraic $k$-group contains a unique maximal multiplicative type subgroup. We show that the converse also holds. More precisely, we prove the following:

\begin{theorem}\label{uniquemaxtoroid} Let $G$ be a connected affine algebraic $k$-group. The following are equivalent:

\vspace{1mm}\noindent (i) $G$ is nilpotent.

\vspace{1mm}\noindent (ii) $G$ contains a unique maximal multiplicative type subgroup.

\vspace{1mm}\noindent (iii) $G$ contains a unique maximal toroid.

\vspace{1mm}\noindent (iv) All toroids of $G$ are central.
\end{theorem} 

\noindent One way of interpreting Theorem \ref{uniquemaxtoroid} is as follows. Consider the subgroup $G_t$ of $G$ which is generated by its toroids. Then Theorem \ref{uniquemaxtoroid} is a first step towards understanding $G_t$: it characterises those algebraic groups $G$ for which $G_t$ is itself a toroid. At the other extreme we have those algebraic groups $G$ which satisfy $G=G_t$, let us call them \textit{generated by toroids}. We would like to characterise such groups. For instance, does $G=G_t$ if $G$ is perfect? How about if $G$ has no non-trivial unipotent subgroups? And finally, is $G_t$ normal in $G$? At present these questions are open, as far as we are aware.

\section{Preliminaries}\label{preliminaries}

\noindent Let $k$ be a field. Let $p\geq 0$ be the characteristic of $k$.

\vspace{2mm}\noindent By an \textit{algebraic $k$-group} we mean a group scheme of finite type over $k$. We do not assume that algebraic $k$-groups are smooth, nor affine (however, in this paper we are only concerned with affine algebraic $k$-groups). By a \textit{subgroup} $H$ of an algebraic $k$-group $G$ we mean a locally closed $k$-subgroup scheme; note that $H$ is automatically closed by \cite[047T]{Stacks} and is of finite type by \cite[01T5, 01T3]{Stacks}. The \textit{image} of a homomorphism of algebraic groups refers to the schematic image, in the sense of \cite[01R7]{Stacks}. For the case where $p>0$, we denote by $\alpha_p$ and $\mu_p$ the first Frobenius kernels of the additive group $\mathbb{G}_a$ and the multiplicative group $\mathbb{G}_m$ respectively.

\vspace{2mm}\noindent Our goal in the remainder of this section is to prove Corollary \ref{unionmulttypetoroidcor}, which is a technical result that we will need to prove Theorem \ref{maxtoroids}.

\vspace{2mm}\noindent First, let us recall the definition of the schematic union. Let $X$ be a scheme. Let $\mathcal{Y}$ be a set of closed subschemes of $X$. [N.B. No set-theoretic issues arise here, as the category of schemes is well-powered.] The \textit{schematic union} $\bigcup_{Y\in\mathcal{Y}}Y$ of $\mathcal{Y}$ in $X$ is defined as follows. The category of schemes admits all -- even uncountable -- coproducts, so we may form the coproduct scheme $\coprod_{Y\in\mathcal{Y}} Y$. Let $\xi:\coprod_{Y\in\mathcal{Y}} Y \to X$ be the morphism of schemes given by inclusion on each component. We then define $\bigcup_{Y\in\mathcal{Y}}Y$ to be the (schematic) image of $\xi$. If $X$ is affine then $\bigcup_{Y\in\mathcal{Y}}Y$ is the closed subscheme of $X$ whose associated ideal is the intersection of those of every $Y\in\mathcal{Y}$.

\vspace{2mm}\noindent Henceforth, let $G$ be an affine algebraic $k$-group. Let $\{H_i \hspace{0.5mm}|\hspace{0.5mm} i \in I\}$ be a -- possibly uncountable -- set of (closed) subgroups of $G$ that is upward directed under inclusion (i.e. $I$ is an upward directed poset where, for $i,j\in I$, $i\leq j$ if and only if $H_i\subseteq H_j$). Let $H:=\bigcup_{i\in I}H_i$, the schematic union of the $H_i$'s in $G$. It was shown in \cite[Prop.\ 2.2]{Se} that this schematic union $H$ is a subgroup of $G$. 

\begin{proposition}\label{commpresschunion} Suppose $H_i$ is commutative for every $i\in I$. Then $H$ is commutative.
\begin{proof} We first consider the case where $I$ is a (possibly uncountable) chain. Taking the centraliser in $G$ of every element of $\{H_i \hspace{0.5mm}|\hspace{0.5mm}i\in I\}$ gives us a downward directed chain \begin{equation}\label{centchain} \big(Z_G(H_i)\big)_{i\in I}.\end{equation} Since $G$ is a Noetherian scheme, 
it satisfies the descending chain condition on closed subschemes. So this chain $(\ref{centchain})$ stabilises. 
That is, there exists a sufficiently large $r\in I$ such that $Z_G(H_i)=Z_G(H_r)$ for all $i\in I$ with $i \geq r$. Consequently, $Z_G(H)=Z_G(H_r)$. 
It follows that \begin{equation}\label{centers} Z(H_i)=H_i\cap Z(H)\end{equation} for every $i\geq r$. 

\vspace{2mm}\noindent Consider the natural projection $\rho:H\to H/Z(H)$. By definition of the schematic union, observe that \begin{equation}\label{imunion}\rho(H)=\mathsmaller{\bigcup}_{i\in I}\rho(H_i).\end{equation} 
Using Equation $(\ref{centers})$, for every $i\geq r$ we have \begin{equation}\label{rhoHi}\rho(H_i)
\cong H_i/(H_i \cap Z(H))=H_i/Z(H_i).\end{equation} By assumption, $H_i$ is commutative for every $i\in I$. Then, by Equations $(\ref{imunion})$ and $(\ref{rhoHi})$, $H$ is also commutative. 

\vspace{2mm}\noindent It remains to consider the general case; that is, $I$ is an arbitrary upward directed poset. Let $\mathcal{Z}$ be the poset of commutative subgroups of $G$, under inclusion. We already showed that every chain in $\mathcal{Z}$ admits a supremum; namely, its schematic union in $G$. It is a basic set-theoretic fact that, given a poset in which every chain admits a supremum, then every upward directed subset admits a supremum. 
Hence this is true of $\mathcal{Z}$. In particular, our upward directed subset $\{H_i \hspace{0.5mm}|\hspace{0.5mm}i\in I\}$ of $\mathcal{Z}$ admits a supremum in $\mathcal{Z}$; let us denote it by $\overline{H}$. By definition of the schematic union, $H$ is contained in $\overline{H}$. But we already know that $H$ is a subgroup of $G$, so $H$ is commutative (and thus $H=\overline{H}$).
\end{proof}
\end{proposition}

\begin{corollary}\label{unionmulttypetoroidcor} Suppose $H_i$ is of multiplicative type (resp. a toroid) for every $i\in I$. Then $H$ is of multiplicative type (resp. a toroid).
\begin{proof} Suppose each $H_i$ is of multiplicative type. Then $H$ is commutative, by Proposition \ref{commpresschunion}. According to \cite[IV, \S3, Thm.\ 1.1(a)]{DG} there exists a largest multiplicative type subgroup of $H$, which we denote by $H_s$. In particular, every $H_i$ is contained in $H_s$. But by definition $H$ is the smallest closed subscheme of $G$ which contains all of the $H_i$'s, so $H=H_s$. That is, $H$ is of multiplicative type. 

\vspace{2mm}\noindent Suppose in addition that each $H_i$ is connected. A similar argument as above shows that $H$ is connected, since there exists a largest connected subgroup of $H$; namely, the connected component of the identity $H^{\circ}$. In other words, $H$ is a toroid.
\end{proof}
\end{corollary}

\section{Maximal toroids and Cartan subgroups}\label{tc}

\noindent In this section we prove Theorems \ref{maxtoroids} and \ref{Cartans}, along with Corollaries \ref{Cartanscor1}, \ref{Cartanscor2}, \ref{Cartanscor3}, \ref{Cartanscor4} and \ref{Cartanscor5}.

\vspace{2mm}\noindent Let $k$ be a field of characteristic $p\geq 0$. Let $\overline{k}$ be an algebraic closure of $k$, and let $k^{\sep}$ be the separable closure of $k$ in $\overline{k}$.

\vspace{2mm}\noindent \underline{Proof of Theorem \ref{maxtoroids}.}

\vspace{2mm}\noindent Let $G$ be an affine algebraic $k$-group.

\begin{proof} (a). Let $S$ be a toroid of $G$. Let $\mathcal{S}$ be the poset of toroids of $G$ which contain $S$, under inclusion. According to Corollary \ref{unionmulttypetoroidcor} every chain in $\mathcal{S}$ has an upper bound in $\mathcal{S}$; namely, its schematic union in $G$. Then applying Zorn's lemma to $\mathcal{S}$, we deduce that $\mathcal{S}$ has a maximal element. The second assertion follows immediately, as the trivial subgroup of $G$ is a toroid.

\vspace{2mm}\noindent (b). Suppose $G$ is smooth. Since $G^{\circ}$ contains all connected subgroups of $G$, without loss of generality we may assume that $G$ is connected. Let $T$ be a maximal toroid of $G$. First consider the case where $T$ is central in $G$. Let $S$ be a maximal torus of $G$. Since $S$ and $T$ commute, $ST$ is again a toroid of $G$. Hence $S\subseteq T$ by maximality. In particular, $G=Z_G(S)$. Then applying \cite[Prop.\ 17.44]{Mi} tells us that $G$ is nilpotent, since both $G$ and $S$ are smooth. Hence $S=T$, by \cite[Thm.\ 16.47(a),(c)]{Mi}.

\vspace{2mm}\noindent Henceforth assume that $T$ is non-central in $G$. We induct on the dimension $\dim G$ of $G$. Assume that $\dim G=1$. Then $G$ is either unipotent or a torus. If $G$ is unipotent then $T$ is trivial, which is a contradiction, so $G=T$. In particular, $T$ is smooth.

\vspace{2mm}\noindent Next assume that $\dim G>1$. Since $T$ is non-central in $G$, its centraliser $Z_G(T)$ is a proper subgroup of $G$. Then $\dim Z_G(T) \lneq \dim G$, as $G$ is smooth and connected. Since $G$ is smooth, applying \cite[XI, Cor.\ 2.4]{SGA3} 
tells us that $Z_G(T)$ (and hence $Z_G(T)^{\circ}$) is smooth. Note that $T$ is a maximal toroid of $Z_G(T)^{\circ}$. Then applying the inductive hypothesis to $Z_G(T)^{\circ}$ tells us that $T$ is smooth.
\end{proof}

\noindent We now move on to the proof of Theorem \ref{Cartans}. We will need the following characterisation of unipotency.

\begin{lemma}\label{nomulttypesubgps} Assume $p>0$. Let $G$ be a connected affine algebraic $k$-group. The following are equivalent:

\vspace{1mm}\noindent (i) $G$ is unipotent.

\vspace{1mm}\noindent (ii) $G_{k^{\sep}}$ does not contain a copy of $\mu_p$.

\vspace{1mm}\noindent (iii) $G$ contains no non-trivial toroids.

\begin{proof} (i) $\implies$ (ii). This is immediate, as the property of unipotency is preserved by base change and by taking subgroups. 

\vspace{2mm}\noindent (ii) $\implies$ (iii). We prove the contrapositive. Suppose there exists a non-trivial toroid $T$ of $G$. Recall from \cite[IV, \S 1, 1.2]{DG} that every height 1 diagonalisable algebraic $k$-group is isomorphic to a product of finitely many copies of $\mu_p$. Hence, since $(T_1)_{k^{\sep}}$ is non-trivial, it contains a copy of $\mu_p$.

\vspace{2mm}\noindent (iii) $\implies$ (i). Suppose $G$ contains no non-trivial toroids. Let $x$ be any non-zero element of the Lie algebra $\Lie(G)$ of $G$. Consider the restricted subalgebra $\langle x \rangle$ of $\Lie(G)$ which is generated by $x$. It is easy to see that $\langle x \rangle$ is commutative. 

\vspace{2mm}\noindent Let $Z$ be the height 1 subgroup of $G$ associated to $\langle x \rangle$ via the Demazure-Gabriel correspondence \cite[II, \S 7, Prop.\ 4.1]{DG}. Note that $Z$ is commutative, as this is true of $\langle x \rangle$. By assumption the largest toroid $Z_s$ of $Z$ is trivial. Hence $Z$ is unipotent, by \cite[IV, \S 3, Thm.\ 1.1(a)]{DG}. This implies that $x$ is $p$-nilpotent, by \cite[XVII, Cor.\ 3.7]{SGA3}. Since $x$ was chosen arbitrarily, we deduce that $\Lie(G)$ is $p$-nilpotent. 
Then applying \cite[XVII, Prop.\ 4.3.1]{SGA3} tells us that $G$ is unipotent. 
\end{proof}
\end{lemma}

\noindent We need two more lemmas.
 
\begin{lemma}\label{Cartansaremax} Let $G$ be an affine algebraic $k$-group. Let $C$ be a Cartan subgroup of $G$. Then $C$ is maximal amongst connected nilpotent subgroups of $G$. 
\begin{proof} Let $H$ be a connected nilpotent subgroup of $G$ which contains $C$. Consider the chain of successive normalisers \begin{equation}\label{sucnorm} N_H(C) \subseteq N_H(N_H(C)) \subseteq...\end{equation} of $C$ in $H$. For every integer $i\geq 1$, we observe that the $i$'th term of this chain $(\ref{sucnorm})$ contains the $i$'th center of $H$ (the argument is the same as for abstract groups). 
Hence, since $H$ is nilpotent, this chain $(\ref{sucnorm})$ terminates after finitely many terms at $H$. 

\vspace{2mm}\noindent Now observe that $N_H(N_H(C))\subseteq N_H(N_H(C)^{\circ})$, since $N_H(C)^{\circ}$ is a characteristic subgroup of $N_H(C)$. But $C$ is a Cartan subgroup of $H$, so $N_H(C)=N_H(N_H(C))$. Consequently, $H=N_H(C)=C$.
\end{proof}
\end{lemma}

\begin{lemma}\label{normchainterm} Let $G$ be an affine algebraic $k$-group. Let $H$ be a subgroup of $G$. There exists a sufficiently large integer $r$ such that $N_G(H_i)=N_G(H^{\circ})$ for every integer $i\geq r$.
\begin{proof} For each integer $i\geq 0$, the $i$'th Frobenius kernel $H_i$ is a characteristic subgroup of $H$. Since $(H_{i+1})_i=H_i$ for every $i\geq 0$, we have a descending chain \begin{equation}\label{normchain} G=N_G(H_0)\supseteq N_G(H_1) \supseteq...\end{equation} of subgroups of $G$. Since $G$ is a Noetherian scheme, this chain $(\ref{normchain})$ stabilises. That is, there exists an integer $r$ such that $N_G(H_r)$ normalises $H_i$ for all $i \geq 0$. It follows from \cite[I, \S 7.17]{Ja} (the density of Frobenius kernels) that $H^{\circ}=\bigcup_{i \geq 0}H_i$, as a schematic union. 
Thus $N_G(H_r)$ normalises $H^{\circ}$. Hence $N_G(H_r)=N_G(H^{\circ})$, since $H_r$ is characteristic in $H^{\circ}$. 
We are done, as $N_G(H_i)=N_G(H_r)$ for all $i\geq r$.
\end{proof}
\end{lemma}

\noindent We now have enough to prove Theorem \ref{Cartans}. 
For the reader's convenience, we restate it here. Recall that $p\geq 0$ is the characteristic of $k$, and $\overline{k}$ is an algebraic closure of $k$.

\vspace{2mm}\noindent \textbf{Theorem 1.2.} Let $G$ be an affine algebraic $k$-group. Let $C$ be a subgroup of $G$. The following are equivalent:

\vspace{1mm}\noindent (i) $C$ is a Cartan subgroup of $G$.

\vspace{1mm}\noindent (ii) $C_{\overline{k}}$ is a Cartan subgroup of $G_{\overline{k}}$.

\vspace{1mm}\noindent (iii) $C$ is connected, and 

\vspace{-1mm}\begin{itemize} \item if $p=0$ the Lie algebra $\Lie(C)$ is a Cartan subalgebra of $\Lie(G)$, \item if $p>0$ there exists a positive integer $r$ such that, for every $i\geq r$, the $i$'th Frobenius kernel $C_i$ is a Cartan subgroup of $G_i$. 
\end{itemize} 

\vspace{1mm}\noindent (iv) $C=Z_G(T)^{\circ}$, for some maximal toroid $T$ of $G$.

\vspace{1mm}\noindent (v) $C$ is a connected nilpotent subgroup of $G$ which is maximal amongst all such subgroups of $G$, and every toroid of $G$ which normalises $C$ must centralise $C$.

\begin{proof} \vspace{2mm}\noindent (i)$\iff$(ii). The property of nilpotency, and the formation of normalisers and the connected component of the identity, are all invariant under base change by algebraic field extensions. So this equivalence is clear.

\vspace{2mm}\noindent (i)$\implies$(iii). Suppose $C$ is a Cartan subgroup of $G$.

\vspace{2mm}\noindent Assume $p=0$. Then $\Lie(C)$ is nilpotent, by \cite[IV, \S 4, Cor.\ 1.6]{DG}. Using \cite[II, \S 6, Prop.\ 2.1(d)]{DG}, we see that $$\Lie(C)=\Lie(N_G(C)^{\circ})=N_{\Lie(G)}(\Lie(C)).$$ That is, $\Lie(C)$ is a Cartan subalgebra of $\Lie(G)$.

\vspace{2mm}\noindent Assume $p>0$. By Lemma \ref{normchainterm}, there exists an integer $r$ such that $N_G(C_i)=N_G(C)$ for every integer $i\geq r$. Hence, for every $i\geq r$, we have $$N_{G_i}(C_i)^{\circ}=N_G(C_i)\cap G_i=
N_G(C)^{\circ}\cap G_i=C\cap G_i=C_i.$$ Certainly $C_i$ is connected and nilpotent for each $i\geq r$, so it is a Cartan subgroup of $G_i$.

\vspace{2mm}\noindent (iv)$\implies$(i). Suppose $C=Z_G(T)^{\circ}$, where $T$ is a maximal toroid of $G$. Let $T'$ be another toroid of $C$. Since $T$ and $T'$ commute, $TT'$ 
is again a toroid of $G$. Hence $T' \subseteq T$ by maximality. That is, all toroids of $C$ are contained in $T$. 

\vspace{2mm}\noindent We claim that the quotient $Q:=C/T$ contains no non-trivial toroids. To see this, let $Y$ be a non-trivial toroid of $Q$. The preimage of $Y$ under the natural projection $\rho:C \to Q$ is also a toroid by \cite[Thm.\ 15.39]{Mi}. But then $\rho^{-1}(Y) \subseteq T$, by the uniqueness of $T$. This is a contradiction, proving the claim. Hence $Q$ is unipotent, by Lemma \ref{nomulttypesubgps}. Therefore $C$ is nilpotent, as it is a central extension of $Q$ by $T$.

\vspace{2mm}\noindent Consider the largest central toroid $Z(C)^{\circ}_s$ of $C$. Note that $T=Z(C)^{\circ}_s$, by maximality of $T$. It follows from \cite[IV, \S 3, Thm.\ 1.1(a)]{DG} that $T$ is a characteristic subgroup of $Z(C)$. Hence $T$ is a characteristic subgroup of $C$, by transitivity. 
Then applying \cite[Cor.\ 12.38]{Mi} to $N_G(C)^{\circ}$ tells us that $T$ is central in $N_G(C)^{\circ}$. Therefore $C=N_G(C)^{\circ}$. 
In summary, we have shown that $C$ is a Cartan subgroup of $G$.

\vspace{2mm}\noindent (iii)$\implies$(i). Suppose (iii) holds. Assume $p=0$. Then $C$ is nilpotent by \cite[IV, \S 4, Cor.\ 1.6]{DG}, and $$\Lie(N_G(C))=N_{\Lie(G)}(\Lie(C))=\Lie(C).$$ by \cite[II, \S 6, Prop.\ 2.1(d)]{DG}. Hence $N_G(C)^{\circ}=C$ by \cite[II, \S 6, Prop.\ 2.1(a)]{DG}.

\vspace{2mm}\noindent Now assume $p>0$. By assumption, there exists an integer $r\geq 1$ such that $C_i$ is a Cartan subgroup of $G_i$ for every $i\geq r$. 
Since each $C_i$ is infinitesimal and nilpotent, by \cite[IV, \S4, Cor.\ 1.11]{DG} we have a chain of toroids \begin{equation}\label{toroidchain} Z(C_1)_s \subseteq Z(C_2)_s \subseteq ...\end{equation} where, for each $i\geq 1$, $Z(C_i)_s$ is the unique maximal toroid of $C_i$. Consider the schematic union $\bigcup_{i\geq 1}Z(C_i)_s$ of this chain $(\ref{toroidchain})$. By Corollary \ref{unionmulttypetoroidcor}, $\bigcup_{i\geq 1}Z(C_i)_s$ is a toroid.

\vspace{2mm}\noindent We claim that $\bigcup_{i\geq 1}Z(C_i)_s$ is a maximal toroid of $C$. To see this, let $S$ be a maximal toroid of $C$ which contains $\bigcup_{i\geq 1}Z(C_i)_s$. Then $S_i=Z(C_i)_s$ for all $i\geq 1$, by maximality of each $Z(C_i)_s$. But then $$S=\mathsmaller{\bigcup}_{i\geq 1}S_i=\mathsmaller{\bigcup}_{i\geq 1}Z(C_i)_s,$$ by \cite[I, \S 7.17]{Ja} (the density of Frobenius kernels). This proves the claim.

\vspace{2mm}\noindent Applying Equation $(\ref{centers})$ (enlarging $r$ if necessary) to the chain of Frobenius kernels $$C_1\subseteq C_2\subseteq ...$$ ensures that $Z(C_i)=Z(C)_i$ for all $i\geq r$. 
Therefore $\bigcup_{i\geq 1}Z(C_i)_s$ is central in $C$, as this is true of each $Z(C_i)_s$. Hence $C$ is nilpotent, by the (already proved) implication (iv)$\implies$(i).

\vspace{2mm}\noindent It remains to check the normaliser condition. By Lemma \ref{normchainterm} (once again enlarging $r$ if necessary), we can assume that $N_G(C_i)=N_G(C)$ for every $i\geq r$. Therefore $$N_G(C)^{\circ}=\mathsmaller{\bigcup}_{i\geq r}N_G(C)_i=\mathsmaller{\bigcup}_{i\geq r}N_{G_i}(C_i)=\mathsmaller{\bigcup}_{i\geq r} C_i=C.$$

\vspace{1.5mm}\noindent (i)$\implies$(iv). Suppose $C$ is a Cartan subgroup of $G$. We claim that $Z(C)^{\circ}_s$ is a maximal toroid of $G$. For the moment, assume the claim holds. Then $Z_G(Z(C)^{\circ}_s)^{\circ}$ is a Cartan subgroup of $G$, by the (already proved) implication (iv)$\implies$(i). But then Lemma \ref{Cartansaremax} implies that $C=Z_G(Z(C)^{\circ}_s)^{\circ}$, and we are done. It remains to prove the claim.

\vspace{2mm}\noindent First assume that $p=0$. Then $G$ and all of its subgroups are smooth, by Cartier's theorem \cite[Thm.\ 3.23]{Mi}. It then follows from \cite[Thm.\ 12.6]{Bo} that $Z(C)^{\circ}_s$ is a maximal torus -- and hence a maximal toroid -- of $G$ (as stated this result is for $k=\overline{k}$, but the proof works over any perfect field $k$). Henceforth assume that $p>0$.

\vspace{2mm}\noindent Next, consider the case where $G$ is infinitesimal and $C$ is unipotent. Note that all subgroups of $G$ are connected. We will show that $C=G$. We proceed by induction on the order $o(G)$ of $G$ (our inductive hypothesis is that $G$ admits no proper unipotent Cartan subgroups). If $o(G)=o(C)$ then of course $C=G$, so $G$ admits no proper unipotent Cartan subgroups by Lemma \ref{Cartansaremax}. Now assume that $o(G)>o(C)$.

\vspace{2mm}\noindent According to \cite[Thm.\ 1.1]{Se}, there exists a largest unipotent normal subgroup $\Rad_u(G)$ of $G$ (this result is not difficult when $G$ is finite). Observe that $\Rad_u(G)\cdot C$ is a unipotent subgroup of $G$, since the property of unipotency is preserved by taking quotients and group extensions. 
Hence $C$ contains $\Rad_u(G)$, again by Lemma \ref{Cartansaremax}. Consider the natural projection $$\pi:G\to G/\Rad_u(G)=:\overline{G}.$$ Denote $\smash{\overline{C}}:=\pi(C)$. Certainly $\smash{\overline{C}}$ is unipotent. Since $C$ contains $\ker\pi=\Rad_u(G)$, it follows that $\smash{\pi(N_G(C))=N_{\overline{G}}(\overline{C})}$. 
Consequently, $\overline{C}$ is a Cartan subgroup of $\overline{G}$. If $\Rad_u(G)$ is non-trivial then $o(\overline{G})\lneq o(G)$, so $\overline{C}=\overline{G}$ by the inductive hypothesis. Hence $C=G$. So we may assume that $\Rad_u(G)$ is trivial.

\vspace{2mm}\noindent Consider the Frobenius kernel $C_1$ of $C$. If $C_1$ is trivial then so is $C$, in which case $C=N_G(C)=G$. So we may assume $C_1$ is non-trivial. Then $N_G(C_1)$ is a proper subgroup of $G$, since $C_1$ is unipotent and $\Rad_u(G)$ is trivial. In particular, $o(N_G(C_1))\lneq o(G)$. Note that $C$ is a unipotent Cartan subgroup of $N_G(C_1)$. 
Hence $C=N_G(C_1)$, by the inductive hypothesis. Then we have $$C_1=C\cap G_1=N_G(C_1)\cap G_1=N_{G_1}(C_1).$$ That is, $C_1$ is a Cartan subgroup of $G_1$.

\vspace{2mm}\noindent Let $\mathfrak{g}$ be the Lie algebra of $G$, and let $\mathfrak{c}$ be the Lie algebra of $C$. Using the Demazure-Gabriel correspondence \cite[II, \S 7, Cor.\ 4.3(d)]{DG}, we see that $\mathfrak{c}$ is a self-normalising subalgebra of $\mathfrak{g}$. Note that $\mathfrak{c}$ is $p$-nilpotent (and hence nilpotent), as $C$ is unipotent. In other words, $\mathfrak{c}$ is a $p$-nilpotent Cartan subalgebra of $\mathfrak{g}$. Then applying \cite[Ch.\ 2, Thm.\ 4.1]{FS} tells us that $\mathfrak{c}=\mathfrak{g}$. Hence $G$ is unipotent, by \cite[XVII, Prop.\ 4.3.1]{SGA3}. Therefore $C=G$, once again using Lemma \ref{Cartansaremax}.

\vspace{2mm}\noindent It remains to consider the general case (for $p>0$). That is, $G$ is an arbitrary affine algebraic $k$-group and $C$ is a Cartan subgroup of $G$. By the (already proved) implication (i)$\implies$(iii), there exists an integer $r$ such that $C_i$ is a Cartan subgroup of $G_i$ for all $i\geq r$.

\vspace{2mm}\noindent Consider the subgroup $\smash{L:=Z_G(Z(C_r)_s)^{\circ}}$ of $G$. Observe that $C_r/Z(C_r)_s$ is a Cartan subgroup of $L_r/Z(C_r)_s$. 
Since $C_r$ is connected and nilpotent, $C_r/Z(C_r)_s$ is unipotent by \cite[IV, \S 4, Cor.\ 1.11]{DG}. In summary, we have found a unipotent Cartan subgroup $C_r/Z(C_r)_s$ of an infinitesimal $k$-group $L_r/Z(C_r)_s$. But we previously showed that an infinitesimal $k$-group admits no proper unipotent Cartan subgroups, hence 
$C_r=L_r$. An identical argument shows that $C_i=L_i$ for all integers $i\geq r$. Then, again using \cite[I, \S 7.17]{Ja} (the density of Frobenius kernels), we deduce that $C=L$. 

\vspace{2mm}\noindent Applying Equation $(\ref{centers})$ to the schematic union of Frobenius kernels $C=\bigcup_{i\geq 0}C_i$ (enlarging $r$ if necessary) tells us that $Z(C)_r=Z(C_r)$. 
Then we have $$(Z(C)^{\circ}_s)_r=Z(C)^{\circ}_s\cap G_r=(Z(C)\cap G_r)_s=Z(C_r)_s.$$ 
Hence $$C\subseteq Z_G(Z(C)^{\circ}_s)^{\circ} \subseteq Z_G(Z(C_r)_s)^{\circ}=L=C.$$ But again using \cite[IV, \S 4, Cor.\ 1.11]{DG}, we know that $Z(C)^{\circ}_s$ is the unique maximal toroid of $C$. Consequently, $Z(C)^{\circ}_s$ is a maximal toroid of $G$.

\vspace{2mm}\noindent (i)$\implies$(v). Suppose $C$ is a Cartan subgroup of $G$. We already showed in Lemma \ref{Cartansaremax} that $C$ is maximal amongst connected nilpotent subgroups of $G$. Now let $S$ be a toroid of $G$ which normalises $C$. Then $S\subseteq N_G(C)^{\circ}=C$, and hence $S$ is central in $C$ by \cite[IV, \S 4, Cor.\ 1.11]{DG}.

\vspace{2mm}\noindent (v)$\implies$(i). Suppose (v) holds. Let $S$ be a toroid of $N_G(C)$. By Theorem \ref{maxtoroids}(a), $S$ is contained in a maximal toroid $T$ of $N_G(C)$. Then $Z_{N_G(C)}(T)^{\circ}$ is a Cartan subgroup of $N_G(C)$, by the implication (iv)$\implies$(i). By assumption $T$ centralises $C$, and $C$ is maximal amongst connected nilpotent subgroups of $G$. It follows that $C=Z_{N_G(C)}(T)^{\circ}$. In particular, $T\subseteq C$. Then $T=Z(C)^{\circ}_s$, by maximality of $T$. Hence $T$ is a characteristic subgroup of $C$, by \cite[IV, \S 4, Cor.\ 1.11]{DG}. Then $T$ is a normal subgroup of $N_G(C)$, and thus a central subgroup of $N_G(C)^{\circ}$ by \cite[Cor.\ 12.38]{Mi}. Consequently, $$N_G(C)^{\circ}=Z_{N_G(C)}(T)^{\circ}=C.$$ That is, $C$ is a Cartan subgroup of $G$.
\end{proof}

\noindent For the subsequent proofs, we continue to assume that $G$ is an affine algebraic $k$-group. 

\vspace{2mm}\noindent \underline{Proof of Corollary \ref{Cartanscor1}.}

\begin{proof} Let $T$ be a maximal toroid of $G$. By (iv)$\implies$(i) of Theorem \ref{Cartans}, $Z_G(T)^{\circ}$ is a Cartan subgroup of $G$. Consider the largest central toroid $Z(Z_G(T)^{\circ})^{\circ}_s$ of $Z_G(T)^{\circ}$. It follows from (i)$\implies$(iv) of Theorem \ref{Cartans} that $Z(Z_G(T)^{\circ})^{\circ}_s$ is a maximal toroid of $G$. But $T \subseteq Z(Z_G(T)^{\circ})^{\circ}_s$, hence $T=Z(Z_G(T)^{\circ})^{\circ}_s$ by maximality of $T$.

\vspace{2mm}\noindent Given a Cartan subgroup $C$ of $G$, a similar argument shows that $Z(C)^{\circ}_s$ is a maximal toroid of $G$ and that $C=Z_G(Z(C)^{\circ}_s)^{\circ}$. This completes the proof.
\end{proof}

\noindent \underline{Proof of Corollary \ref{Cartanscor2}.}

\begin{proof} By Theorem \ref{maxtoroids}(a), $G$ contains at least one maximal toroid $T$. Then $Z_G(T)^{\circ}$ is a Cartan subgroup of $G$, by Corollary \ref{Cartanscor1}.
\end{proof}

\noindent \underline{Proof of Corollary \ref{Cartanscor3}.}

\begin{proof} Suppose $G$ is smooth. Let $C$ be a Cartan subgroup of $G$. Observe that $C=Z_G(Z(C)^{\circ}_s)^{\circ}$, by Corollary \ref{Cartanscor1}. Then applying \cite[XI, Cor.\ 2.4]{SGA3} tells us that $C$ is smooth, since $G$ is smooth and $Z(C)^{\circ}_s$ is of multiplicative type.
\end{proof}

\noindent \underline{Proof of Corollary \ref{Cartanscor4}.}

\begin{proof} Suppose $T$ is a maximal toroid of $G$. By Corollary \ref{Cartanscor1}, $Z_G(T)^{\circ}$ is a Cartan subgroup of $G$ and $T=Z(Z_G(T)^{\circ})^{\circ}_s$. Then $\smash{Z_G(T)^{\circ}_{\overline{k}}}$ is a Cartan subgroup of $G_{\overline{k}}$, by (i)$\implies$(ii) of Theorem \ref{Cartans}. Note that the formation of the largest central toroid commutes with base change by $\overline{k}/k$, hence $\smash{T_{\overline{k}}=Z(Z_G(T)^{\circ}_{\overline{k}})^{\circ}_s}$. Thus $T_{\overline{k}}$ is a maximal toroid of $G_{\overline{k}}$, again using Corollary \ref{Cartanscor1}. The converse is obvious.
\end{proof}

\noindent \underline{Proof of Corollary \ref{Cartanscor5}.}

\begin{proof} Suppose $T$ is a maximal toroid of $G$. By Corollary \ref{Cartanscor1}, $Z_G(T)^{\circ}$ is a Cartan subgroup of $G$ and $T=Z(Z_G(T)^{\circ})^{\circ}_s$. Since $Z_G(T)^{\circ}$ is nilpotent, by \cite[IV, \S 4, Cor.\ 1.11]{DG} all toroids of $Z_G(T)^{\circ}$ are central and hence are contained in $T$. 
Consequently, for every integer $i\geq 1$ the infinitesimal toroid $Z(Z_G(T)_i)_s$ is contained in $T\cap G_i=T_i$. By (i)$\implies$(iii) of Theorem \ref{Cartans}, there exists a positive integer $r$ such that $Z_G(T)_i$ is a Cartan subgroup of $G_i$ for all $i\geq r$. Now fix an integer $i\geq r$. Then $Z(Z_G(T)_i)_s$ is a maximal toroid of $G_i$, again by Corollary \ref{Cartanscor1}. But we showed that $Z(Z_G(T)_i)_s\subseteq T_i$, hence $Z(Z_G(T)_i)_s=T_i$ by maximality.

\vspace{2mm}\noindent For the converse, suppose $T$ is a connected subgroup of $G$, and that there exists a positive integer $r$ such that $T_i$ is a maximal toroid of $G_i$ for every $i\geq r$. Then of course $T_1$ is a toroid, hence $T$ is a toroid by \cite[Cor.\ 1.6]{Se}. Let $S$ be a toroid of $G$ which contains $T$. Then $S_i=T_i$ for all $i\geq r$, by maximality of each $T_i$. But then $$S=\mathsmaller{\bigcup}_{i\geq r}S_i=\mathsmaller{\bigcup}_{i\geq r}T_i=T,$$ by the density of Frobenius kernels \cite[I, \S 7.17]{Ja}.
\end{proof}

\section{Generation problems}\label{gc}

\noindent In this section we prove Theorems \ref{Gcsubgroup} and \ref{uniquemaxtoroid}.

\vspace{2mm}\noindent Recall that $k$ is a field of characteristic $p\geq 0$. Let $\overline{k}$ be an algebraic closure of $k$.

\begin{lemma}\label{Cartansgrow} Let $G$ be an affine algebraic $k$-group. Let $S$ be a toroid of $G$, and let $C$ be a Cartan subgroup of $Z_G(S)$. Then $C$ is contained in a Cartan subgroup of $G$.
\begin{proof} According to Corollary \ref{Cartanscor1}, $Z(C)^{\circ}_s$ is a maximal toroid of $Z_G(S)$. Observe that $Z(C)^{\circ}_s\cdot S$ is a toroid of $Z_G(S)$, since $S$ is central in $Z_G(S)$. 
Hence $Z(C)^{\circ}_s$ contains $S$, by maximality.

\vspace{2mm}\noindent Now let $T$ be a maximal toroid of $G$ which contains $Z(C)^{\circ}_s$. Again by Corollary \ref{Cartanscor1}, we know that $Z_G(T)^{\circ}$ is a Cartan subgroup of $G$. Note that $T\subseteq Z_G(S)$, since $S\subseteq Z(C)^{\circ}_s\subseteq T$. Hence $T=Z(C)^{\circ}_s$, by maximality. Therefore $C\subseteq Z_G(T)^{\circ}$.
\end{proof}
\end{lemma}

\noindent Henceforth assume $k$ has positive characteristic $p>0$. We turn our attention to restricted Lie algebras over $k$. We say that a restricted Lie algebra over $k$ is \textit{toral} if it is abelian and all of its elements are semisimple. [N.B. Many authors instead use the term \textit{torus}, but we prefer to avoid this term to avoid any confusion with algebraic groups.]

\begin{corollary}\label{Cartansgrowcor} Let $\mathfrak{g}$ be a finite-dimensional restricted Lie algebra over $k$. Let $s$ be a semisimple element of $\mathfrak{g}$. Let $\mathfrak{c}$ be a Cartan subalgebra of the centraliser $\mathfrak{z}_{\mathfrak{g}}(s)$. Then $\mathfrak{c}$ is contained in a Cartan subalgebra of $\mathfrak{g}$.
\begin{proof} Let $G$ be the height 1 algebraic $k$-group associated to $\mathfrak{g}$ under the Demazure-Gabriel correspondence \cite[II, \S 7, no.\ 4]{DG}. Note that toroids (resp. Cartan subgroups) of $G$ correspond to toral subalgebras (resp. Cartan subalgebras) of $\mathfrak{g}$. The (restricted) subalgebra of $\mathfrak{g}$ generated by $s$ is toral, since $s$ is semisimple. The result then follows from combining Lemma \ref{Cartansgrow} with the Demazure-Gabriel correspondence.
\end{proof}
\end{corollary}

\noindent The following key result appears to be known, at least over an infinite field $k$, yet we are unaware of an explicit statement anywhere in the literature. A similar statement -- but with an additional splitness assumption -- was proved in \cite[Cor.\ 3.11]{Wi}. Alexander Premet has privately given us a proof which works over any infinite field $k$, by showing that the set of regular elements of a finite-dimensional Lie algebra over $k$ is Zariski-dense. 
Here we give a different proof, which works over any field of positive characteristic (including finite fields).

\begin{lemma}\label{genbyCartansubalg} Let $\mathfrak{g}$ be a finite-dimensional restricted Lie algebra over $k$. Then $\mathfrak{g}$ is spanned by its Cartan subalgebras.

\begin{proof} Recall from \cite[Ch.\ 2, Thm.\ 4.1]{FS} that the Cartan subalgebras of $\mathfrak{g}$ are precisely the centralisers of maximal toral subalgebras of $\mathfrak{g}$; we will repeatedly use this result throughout this proof.

\vspace{2mm}\noindent We induct on the dimension $\dim\mathfrak{g}$ of $\mathfrak{g}$. If $\dim\mathfrak{g}=1$ then $\mathfrak{g}$ is commutative, so it is a Cartan subalgebra of itself. Henceforth assume $\dim\mathfrak{g}>1$. Let $\mathfrak{g}_c$ denote the subspace of $\mathfrak{g}$ spanned by its Cartan subalgebras.

\vspace{2mm}\noindent We consider the following three cases separately.

\vspace{2mm}\noindent \underline{\textbf{Case 1.}} Suppose the center $\mathfrak{z}(\mathfrak{g})$ of $\mathfrak{g}$ is non-zero. 

\vspace{2mm}\noindent By the inductive hypothesis, $\mathfrak{g}/\mathfrak{z}(\mathfrak{g})$ is spanned by its Cartan subalgebras. Consider the natural projection $$\zeta:\mathfrak{g}\to \mathfrak{g}/\mathfrak{z}(\mathfrak{g}).$$ It is easy to see that the preimage of a Cartan subalgebra under a central quotient map remains a Cartan subalgebra. Therefore $\mathfrak{g}=\mathfrak{g}_c$.

\vspace{2mm}\noindent \underline{\textbf{Case 2.}} Suppose the $p$-nilpotent radical $\Rad_u(\mathfrak{g})$ of $\mathfrak{g}$ is non-zero (i.e. there exists a non-zero $p$-nilpotent ideal of $\mathfrak{g}$). 

\vspace{2mm}\noindent We will need the following fact: given a surjective homomorphism $\rho:\mathfrak{g}\to\mathfrak{g}'$ of finite-dimensional restricted Lie algebras, for every Cartan subalgebra $\mathfrak{c}'$ of $\mathfrak{g}'$ there exists a Cartan subalgebra $\mathfrak{c}$ of $\mathfrak{g}$ such that $\rho(\mathfrak{c})=\mathfrak{c}'$. This fact follows from \cite[Thm.\ 4.4.5.1]{Wi1}. 

\vspace{2mm}\noindent Let $\mathfrak{v}$ be a (non-zero) minimal $p$-nilpotent ideal of $\mathfrak{g}$, and consider the natural projection $$\pi:\mathfrak{g}\to\mathfrak{g}/\mathfrak{v}.$$ By the inductive hypothesis, $\mathfrak{g}/\mathfrak{v}$ is spanned by its Cartan subalgebras. Then applying the aforementioned fact to $\pi$ implies that $\mathfrak{g}$ is spanned by $\mathfrak{g}_c$ along with $\mathfrak{v}$. If $\mathfrak{v}\subseteq\mathfrak{g}_c$ then of course $\mathfrak{g}=\mathfrak{g}_c$, and we are done. So assume (for a contradiction) that $\mathfrak{v}\not\subseteq\mathfrak{g}_c$. We consider the following two possibilities separately: either $\mathfrak{g}/\mathfrak{v}$ is toral, or it isn't.

\vspace{2mm}\noindent First, suppose $\mathfrak{g}/\mathfrak{v}$ is not toral. Let $\mathfrak{t}'$ be a maximal toral subalgebra of $\mathfrak{g}/\mathfrak{v}$. According to \cite[Ch.\ 2, Thm.\ 4.5(2)]{FS}, any maximal toral subalgebra of $\pi^{-1}(\mathfrak{t}')$ is a maximal toral subalgebra of $\mathfrak{g}$. Consequently, any Cartan subalgebra of $\pi^{-1}(\mathfrak{t}')$ is contained in a Cartan subalgebra of $\mathfrak{g}$. By the inductive hypothesis, $\pi^{-1}(\mathfrak{t}')$ is spanned by its Cartan subalgebras. Hence $\mathfrak{v}\subseteq\pi^{-1}(\mathfrak{t}')\subseteq \mathfrak{g}_c$, which is a contradiction. 

\vspace{2mm}\noindent Next, suppose $\mathfrak{g}/\mathfrak{v}$ is toral. Again using \cite[Ch.\ 2, Thm.\ 4.5(2)]{FS}, there exists a maximal toral subalgebra $\mathfrak{t}$ of $\mathfrak{g}$ such that $\mathfrak{g}=\mathfrak{v}\rtimes\mathfrak{t}$. By minimality of $\mathfrak{v}$ we deduce that $\mathfrak{v}$ is abelian and $\mathfrak{v}^p=0$ (one way to see this is to apply \cite[XVII, Lem.\ 3.9(i)]{SGA3} to $\mathfrak{v}$, via the Demazure-Gabriel correspondence \cite[II, \S 7, no.\ 4]{DG}). 
It follows that $\mathfrak{g}_c\cap\mathfrak{v}$ is an ideal of $\mathfrak{g}$, since $\mathfrak{g}$ is spanned by $\mathfrak{g}_c$ and $\mathfrak{v}$. 
Then $\mathfrak{g}_c\cap\mathfrak{v}=0$, since $\mathfrak{v}\not\subseteq\mathfrak{g}_c$ and again by minimality of $\mathfrak{v}$. Hence $\mathfrak{g}_c=\mathfrak{t}$, by dimension considerations. However, it is easy to find a semisimple element of $\mathfrak{g}$ which is not contained in $\mathfrak{t}$. For instance, let $t\in\mathfrak{t}$ and $v\in\mathfrak{v}$ such that $[t,v]\neq 0$. A simple calculation shows us that $t+[t,v]$ is semisimple. 
Therefore $t+[t,v]\in\mathfrak{g}_c=\mathfrak{t}$, and hence $[t,v]\in\mathfrak{t}\cap\mathfrak{v}=0$. 
Again, we have a contradiction.

\vspace{2mm}\noindent \underline{\textbf{Case 3.}} Suppose $\mathfrak{z}(\mathfrak{g})=0$ and $\Rad_u(\mathfrak{g})=0$. 

\vspace{2mm}\noindent Let $x\in\mathfrak{g}$. Our goal is to show that $x\in\mathfrak{g}_c$. Consider the centraliser $\mathfrak{z}_\mathfrak{g}(x)$ of $x$ in $\mathfrak{g}$.

\vspace{2mm}\noindent First, suppose $\mathfrak{z}_\mathfrak{g}(x)$ is not $p$-nilpotent. Then by \cite[Ch.\ 2, Thm.\ 3.4]{FS} 
there exists a non-zero semisimple element $s\in\mathfrak{z}_\mathfrak{g}(x)$. Note that $\mathfrak{z}_\mathfrak{g}(s)$ is a proper subalgebra of $\mathfrak{g}$, since $\mathfrak{z}(\mathfrak{g})=0$. Hence $\mathfrak{z}_\mathfrak{g}(s)$ is spanned by its Cartan subalgebras, by the inductive hypothesis. Then applying Corollary \ref{Cartansgrowcor} tells us that $\mathfrak{z}_\mathfrak{g}(s)\subseteq\mathfrak{g}_c$. So $x\in \mathfrak{g}_c$.

\vspace{2mm}\noindent Henceforth suppose $\mathfrak{z}_\mathfrak{g}(x)$ is $p$-nilpotent. Let $\mathfrak{u}$ be a maximal $p$-nilpotent subalgebra of $\mathfrak{g}$ which contains $\mathfrak{z}_\mathfrak{g}(x)$. Note that $\mathfrak{u}$ is not an ideal of $\mathfrak{g}$, since $\Rad_u(\mathfrak{g})=0$. So its normaliser $\mathfrak{n}_{\mathfrak{g}}(\mathfrak{u})$ is a proper subalgebra of $\mathfrak{g}$. Hence, by the inductive hypothesis, $\mathfrak{n}_{\mathfrak{g}}(\mathfrak{u})$ is spanned by its Cartan subalgebras. 

\vspace{2mm}\noindent Let $\mathfrak{c}$ be a Cartan subalgebra of $\mathfrak{n}_{\mathfrak{g}}(\mathfrak{u})$. Assume (for a contradiction) that $\mathfrak{c}$ is $p$-nilpotent. Then $\mathfrak{n}_{\mathfrak{g}}(\mathfrak{u})$ is also $p$-nilpotent, and hence $\mathfrak{u}=\mathfrak{n}_{\mathfrak{g}}(\mathfrak{u})$ by maximality of $\mathfrak{u}$. In other words, $\mathfrak{u}$ is a Cartan subalgebra of $\mathfrak{g}$. But then the same argument tells us that $\mathfrak{g}$ itself is $p$-nilpotent, which contradicts the fact that $\Rad_u(\mathfrak{g})=0$. So $\mathfrak{c}$ is not $p$-nilpotent. Consequently, there exists a non-zero semisimple element $t\in\mathfrak{z}(\mathfrak{c})$. Note that $\mathfrak{z}_{\mathfrak{g}}(t)$ is a proper subalgebra of $\mathfrak{g}$, since $\mathfrak{z}(\mathfrak{g})=0$. Then applying the inductive hypothesis to $\mathfrak{z}_{\mathfrak{g}}(t)$ and using Corollary \ref{Cartansgrowcor}, we deduce that $\mathfrak{c}\subseteq\mathfrak{z}_{\mathfrak{g}}(t)\subseteq\mathfrak{g}_c$. 
Hence $\mathfrak{n}_{\mathfrak{g}}(\mathfrak{u})\subseteq\mathfrak{g}_c$, as it is spanned by its Cartan subalgebras. So indeed $x\in\mathfrak{g}_c$.
\end{proof}
\end{lemma}

\noindent We now have enough to prove Theorems \ref{Gcsubgroup} and \ref{uniquemaxtoroid}.

\vspace{2mm}\noindent \underline{Proof of Theorem \ref{Gcsubgroup}.}

\vspace{2mm}\noindent Assume $p>0$, and let $G$ be an affine algebraic $k$-group. Recall that $G_c$ denotes the subgroup of $G$ which is generated by its Cartan subgroups. 

\begin{proof} Without loss of generality we can and do assume $G$ is connected. Let $\mathfrak{g}$ be the Lie algebra of $G$. First, consider the case where $\mathfrak{g}$ is $p$-nilpotent. Then $G$ is unipotent by \cite[XVII, Prop.\ 4.3.1]{SGA3}. Hence $G=G_c$, and we are done. Henceforth assume $\mathfrak{g}$ is not $p$-nilpotent.

\vspace{2mm}\noindent By the Demazure-Gabriel correspondence \cite[II, \S 7, no.\ 4]{DG}, it suffices to show that $\Lie(G_c)=\mathfrak{g}$. We prove this by induction on the dimension $\dim\mathfrak{g}$ of $\mathfrak{g}$. Assume $\dim\mathfrak{g}=1$. Then $\mathfrak{g}$ is toral since it is not $p$-nilpotent. Hence $G$ is a toroid by \cite[Cor.\ 1.6]{Se}. Therefore $G=G_c$, so of course $\Lie(G_c)=\mathfrak{g}$. Henceforth assume $\dim\mathfrak{g}>1$.

\vspace{2mm}\noindent Next, consider the case where $\mathfrak{g}$ contains no non-zero central semisimple elements. Let $s$ be a non-zero semisimple element of $\mathfrak{g}$ (such an element indeed exists since $\mathfrak{g}$ is not $p$-nilpotent, by \cite[Ch.\ 2, Thm.\ 3.4]{FS}). Consider its centraliser $Z_G(s)$ in $G$. By Lemma \ref{Cartansgrow}, every Cartan subgroup of $Z_G(s)$ is contained in some Cartan subgroup of $G$. Hence $Z_G(s)_c\subseteq G_c$. By assumption $s$ is not central in $\mathfrak{g}$, so $\mathfrak{z}_{\mathfrak{g}}(s)$ is a proper subalgebra of $\mathfrak{g}$. Note that $\Lie(Z_G(s))=\mathfrak{z}_{\mathfrak{g}}(s)$ by \cite[II, \S 7, Cor.\ 4.3(d)]{DG}. Therefore $\Lie(Z_G(s)_c)=\mathfrak{z}_{\mathfrak{g}}(s)$, by the inductive hypothesis. If some Cartan subalgebra $\mathfrak{c}$ of $\mathfrak{g}$ is $p$-nilpotent then $\mathfrak{c}=\mathfrak{g}$ by \cite[Ch.\ 2, Thm.\ 4.1]{FS}, which contradicts our assumption. So none of the Cartan subalgebras of $\mathfrak{g}$ are $p$-nilpotent. We showed in Lemma \ref{genbyCartansubalg} that $\mathfrak{g}$ is spanned by its Cartan subalgebras. It follows that $\mathfrak{g}$ is spanned by the centralisers $\mathfrak{z}_{\mathfrak{g}}(s)$, as $s$ runs over all non-zero semisimple elements of $\mathfrak{g}$. Since our $s$ was chosen arbitrarily, we deduce that $\Lie(G_c)=\mathfrak{g}$.

\vspace{2mm}\noindent It remains to consider the general case. Let $$\zeta:G\to G/Z(G)^{\circ}_s=:G'$$ be the natural projection. Observe that $G'$ contains no non-trivial central toroids, as otherwise we would violate rigidity \cite[Cor.\ 12.38]{Mi}. This implies that $\Lie(G')$ contains no non-zero central semisimple elements. But we already proved this special case, so $G'_c$ contains $G'_1$.

\vspace{2mm}\noindent Given any Cartan subgroup $C'$ of $G'$, observe that $\zeta^{-1}(C')$ is a Cartan subgroup of $G$ (since $\zeta$ is central). It follows that $G_c\supseteq\zeta^{-1}(G'_c)$. Hence $G_c$ contains $\zeta^{-1}(G'_1)$, which in turn contains $G_1$. This completes the proof.
\end{proof}

\noindent \underline{Proof of Theorem \ref{uniquemaxtoroid}.}

\vspace{2mm}\noindent Let $G$ be a connected affine algebraic $k$-group.

\begin{proof} (i) $\implies$ (ii) and (i) $\implies$ (iv). Suppose $G$ is nilpotent. Then by \cite[IV, \S4, Cor.\ 1.11]{DG} there exists a unique maximal multiplicative type subgroup of $G$, namely $Z(G)_s$.

\vspace{2mm}\noindent (ii) $\implies$ (iii). Suppose $G$ contains a unique maximal multiplicative type subgroup $M$. Then certainly $G$ contains a unique maximal toroid, namely the identity component $M^{\circ}$.

\vspace{2mm}\noindent (iv) $\implies$ (iii). Suppose that all toroids of $G$ are central in $G$. Then there exists a unique maximal toroid of $G$, namely $Z(G)^{\circ}_s$.

\vspace{2mm}\noindent (iii) $\implies$ (i). Suppose $G$ contains a unique maximal toroid; let us call it $T$. We consider the cases $p=0$ and $p>0$ separately.

\vspace{2mm}\noindent Assume $p=0$. Then both $G$ and $T$ are smooth, by Cartier's theorem \cite[Thm.\ 3.23]{Mi}. In other words, $T$ is a maximal torus of $G$. Any two maximal tori of $G$ are conjugate by an element of $G(\overline{k})$, by \cite[Cor.\ 11.3(1)]{Bo}. 
Note that $G(\overline{k})$ is Zariski-dense in $G$, as $G$ is smooth. It follows that $T$ is a normal subgroup of $G$. Then $T$ is central in $G$, by \cite[Cor.\ 12.38]{Mi}. That is, $G$ is a Cartan subgroup of itself. In particular, $G$ is nilpotent.

\vspace{2mm}\noindent Henceforth assume $p>0$. We first consider the special case where $G$ contains no non-trivial central toroids. In this case, we claim that $G$ is unipotent.

\vspace{2mm}\noindent Since $G$ has a unique maximal toroid $T$, its Frobenius kernel $T_1$ is the unique maximal toroid of $G_1$. Hence $G_1$ has a unique Cartan subgroup by Corollary \ref{Cartanscor1}, namely $\smash{Z_{G_1}(T_1)}$. Then applying Theorem \ref{Gcsubgroup} to $G_1$ tells us that $\smash{G_1=Z_{G_1}(T_1)}$. Consequently, $G_1$ is nilpotent and satisfies $Z(G_1)_s=T_1$. So $T_1$ is a characteristic subgroup of $Z(G_1)$, by \cite[IV, \S 3, Thm.\ 1.1(a)]{DG}. It follows that $T_1$ is a characteristic subgroup of $G$. Hence $T_1$ is central in $G$, by \cite[Cor.\ 12.38]{Mi}. But by assumption $G$ contains no non-trivial central toroids, so $T_1$ is trivial. Therefore $T$ is trivial, which implies that $G$ is unipotent by Lemma \ref{nomulttypesubgps}. This proves the claim.

\vspace{2mm}\noindent It remains to consider the general case. Let $$\zeta:G\to G/Z(G)^{\circ}_s=:G'$$ be the natural projection. By \cite[Cor.\ 12.42]{Mi} the property of being a toroid is preserved by taking group extensions, and obviously it is preserved by taking quotients. It follows that $\zeta(T)=:T'$ is the unique maximal toroid of $G'$. As already observed in the proof of Theorem \ref{Gcsubgroup}, $G'$ contains no non-trivial central toroids. Hence $G'$ is unipotent by the claim. Therefore $G$ is nilpotent, as it is a central extension of $G'$. This completes the proof of Theorem \ref{uniquemaxtoroid}.
\end{proof}

\noindent An affine algebraic $k$-group $G$ is said to be \textit{generated by Cartans} if $G=G_c$. We conclude this paper by giving an example of a connected $G$ which is \textbf{not} generated by Cartans. 

\begin{example}\label{notgenbyCartans} Assume $p>0$. Let $G=\mathbb{G}_a\rtimes\mu_p$, where the action of $\mu_p$ on $\mathbb{G}_a$ is non-trivial. Let $T$ be a maximal toroid of $G$. Observe that $T\cong\mu_p$, and that $Z_G(T)_1=T$. Then $Z_G(T)^{\circ}$ is a toroid, by \cite[Cor.\ 1.6]{Se}. Hence $T=Z_G(T)^{\circ}$, by maximality. Since $T$ has height 1, the subgroup $G_c$ of $G$ generated by all of its Cartan subgroups equals $G_1$. 
\end{example}

\begin{remark}\label{notgenbyCartansrem} It is known that Cartan subgroups are preserved by surjections between smooth affine $k$-groups, see \cite[11.14, Cor.\ 2]{Bo}. Explicitly, given a surjective homomorphism between smooth affine $k$-groups $f:G\to G'$ and a Cartan subgroup $C$ of $G$, then $f(C)$ is a Cartan subgroup of $G'$. Without the smoothness assumption, Cartan subgroups are not necessarily preserved by faithfully flat homomorphisms. To see this, let $G$ be as in Example \ref{notgenbyCartans}. Consider the second Frobenius kernel $G_2\cong\alpha_{p^2}\rtimes\mu_p$ of $G$ and its derived subgroup $\mathscr{D}(G_2)\cong\alpha_p$. Every Cartan subgroup of $G_2$ is isomorphic to $\mu_p$, so the natural projection $G_2\to G_2/\mathscr{D}(G_2)\cong \alpha_p\times\mu_p$ does not send any Cartan subgroup of $G_2$ onto $G_2/\mathscr{D}(G_2)$ (which is a Cartan subgroup of itself).
\end{remark}

\bigskip\noindent
{\textbf {Acknowledgements}}:
The author would like to thank Wolfgang Soergel, David Stewart, and Alexander Premet for some helpful comments and discussions. In particular, Alexander Premet privately showed us an alternate proof of Lemma \ref{genbyCartansubalg} for infinite fields. The author is employed/funded by the University of Freiburg Mathematical Institute.

\bibliographystyle{amsalpha}

\newcommand{\etalchar}[1]{$^{#1}$}
\providecommand{\bysame}{\leavevmode\hbox to3em{\hrulefill}\thinspace}
\providecommand{\MR}{\relax\ifhmode\unskip\space\fi MR }
\providecommand{\MRhref}[2]{%
	\href{http://www.ams.org/mathscinet-getitem?mr=#1}{#2} }
\providecommand{\href}[2]{#2}

\end{document}